\documentclass{amsart}

\usepackage{amsfonts}
\usepackage{amsmath}
\usepackage{amsthm}
\usepackage{amssymb}
\usepackage[english]{babel}
\usepackage{graphics}
\usepackage{latexsym}
\usepackage{longtable} \usepackage{mathrsfs}
\usepackage{graphicx}
\usepackage{wrapfig} 
\usepackage[pdftex,colorlinks=true]{hyperref}

\newtheorem{theorem}{Theorem}
\newtheorem{corollary}[theorem]{Corollary}
\newtheorem{lemma}[theorem]{Lemma}

\newtheorem{claim}[theorem]{Claim}
\newtheorem{example}[theorem]{Example}

\theoremstyle{definition}
\newtheorem{definition}[theorem]{Definition}
\newtheorem{remark}[theorem]{Remark}

\newcommand{\R}{\mathbb{R}}

\newcommand{\noi}{\noindent}
\newcommand{\ms}{\medskip}

\newcommand{\al}{\alpha}
\newcommand{\be}{\beta}
\newcommand{\ga}{\gamma}

\newcommand{\De}{\Delta}

\newcommand{\Om}{\Omega}

\newcommand{\larrow}{\longrightarrow}
\newcommand{\ot}{\otimes}

\newcommand{\p}{\partial}
\newcommand{\sub}{\subseteq}
\newcommand{\set}{\setminus}
\newcommand{\by}{\times}

\newcommand{\co}{\overline{\textrm{co}}}

\newcommand{\bt}{\begin{theorem}}\newcommand{\et}{\end{theorem}}
\newcommand{\bd}{\begin{definition}}\newcommand{\ed}{\end{definition}}
\newcommand{\bl}{\begin{lemma}}\newcommand{\el}{\end{lemma}}
\newcommand{\beq}{\begin{equation}}\newcommand{\eeq}{\end{equation}}
\newcommand{\bc}{\begin{claim}}\newcommand{\ec}{\end{claim}}
\newcommand{\bex}{\begin{example}}\newcommand{\eex}{\end{example}}
\newcommand{\bcor}{\begin{corollary}}\newcommand{\ecor}{\end{corollary}}
\newcommand{\bp}{\begin{proof}}\newcommand{\ep}{\end{proof}}

\numberwithin{equation}{section}

\begin{document}

\title[Remarks on the Maximum Principle for the $\infty$-Laplacian]{Remarks on the Validity of the Maximum Principle for the $\infty$-Laplacian}

\author{Nikos Katzourakis and Juan Manfredi}
\address{Department of Mathematics and Statistics, University of Reading, Whiteknights, PO Box 220, Reading RG6 6AX, Berkshire, UK} 
\email{n.katzourakis@reading.ac.uk}
\address{Department of Mathematics, University of Pittsburgh, USA}
\email{manfredi@pitt.edu}


\date{}



\begin{abstract} In this note we give three counter-examples which show that the Maximum Principle generally fails for classical solutions of a system and a single equation related to the $\infty$-Laplacian. The first is the tangential part of the $\infty$-Laplace system and the second is the scalar $\infty$-Laplace equation perturbed by a linear gradient term. The interpretations of the Maximum Principle for the system are that of the Convex Hull Property and also of the Maximum Principle of the modulus of the solution.

\end{abstract}

\maketitle

\section{Introduction}

Given a smooth map  $u:\Om \sub \R^n \!\larrow \R^N$, the $\infty$-Laplacian is the PDE system
\beq \label{1}
\De_\infty u \, :=\, \Big(Du \ot Du +  |Du|^2[Du]^\bot \! \ot I\Big):D^2 u \, =  \, 0. 
\eeq
In the above, $u=(u_1,...,u_N)^\top$, $D_i\equiv \frac{\p}{\p x_i}$ and the gradient matrix is understood as 
\[
Du(x)\, =\, (D_iu_\al(x))^{\al =1,...,N}_{i =1,...,n}, 
\]
for $n,N\geq 1$. The system in  index form reads
\[
D_i u_\al\, D_ju_\be \,D^2_{ij} u_\be \ +\ |Du|^2 [Du]_{\al \be}^\bot D^2_{ii} u_\be\ = \ 0,
\]
where triple summation in $i,j \in \{1,...,n\}$ and $\be \in \{1,...,N\}$ is implied although it is not explicitely written. The symbol $|\cdot|$ denotes the Euclidean norm on $\R^{N \by n}$, that is 
\[
|Du|^2\, =\, D_iu_\al D_iu_\al \, =\, Du : Du
\]
and $[Du(x)]^\bot$ denotes the orthogonal projection in $\R^N$ on the nullspace of the transpose of the linear map 
\[
Du(x)\ : \ \R^n \larrow \R^N, \ \ \ \ x\in \Om. 
\]
Geometrically, $[Du(x)]^\bot$ defines the projection on the normal space of the image of the solution at $u(x)$. The hessian of $u$ is viewed as a map $D^2u : \Om\sub \R^n \larrow \R^{Nn^2}$ and
\[
D^2u(x)\, =\, (D^2_{ij}u_\al(x))^{\al =1,...,N}_{i,j =1,...,n}, 
\]
Because of mutual perpendicularity of the two summands comprising the $\infty$-Laplace system, the system actually consists of a pair of two independent systems:
\[
\left\{
\begin{array}{r}
Du\ot Du :D^2u\, =\, 0,\,\ms\\
|Du|^2[Du]^\bot \De u\, =\, 0.
\end{array}
\right.
\]
The $\infty$-Laplace system arises as the analogue of the Euler-Lagrange equation of the supremal functional
\[
E_{\infty}(u,\Om)\, :=\, \|Du\|_{L^\infty(\Om)},
\]
where the $L^\infty$ norm above is interpreted as the essential supremum of the Euclidean norm of the gradient. The scalar case of $\De_\infty$ and $E_\infty$ have a long history and first arose in the papers of Aronsson \cite{A1,A2}. Both the functional and the single equation have been extensively studied ever since, see for instance \cite{C}, \cite{BEJ}, \cite{J}, \cite{CGW} and for a pedagogical introduction see \cite{K}. When $N=1$, the $\infty$-Laplacian simplifies to
\[
Du \ot Du :D^2u\, =\ D_iu\, D_ju\, D^2_{ij}u\, =\, 0
\] 
and the perpendicular term $|Du|^2[Du]^\bot \De u=0$ vanishes identically. The development of the vectorial case of $N\geq 2$ is much more recent and is due to the first author (see \cite{K2}-\cite{K7}). The derivation of the $\infty$-Laplace system from the $L^\infty$ functional and other analytical results have first been established in \cite{K2}. Almost simultaneously to \cite{K2}, Sheffield and Smart studied in \cite{SS} the problem of vectorial Lipschitz extensions which amounts to changing from the Euclidean norm we are using herein to the operator norm on the matrix space of the gradients (this leads to a much more complicated version of ``$\infty$-Laplacian" than ours which may have even multi-valued coefficients).

The vectorial case is largely under development and still not fully understood. Particular difficulties in the study of the system \eqref{1} are the emergence of highly singular solutions and that the operator definining $\De_\infty$ may have discontinuous coefficients even when it is applied to smooth maps (see \cite{K2, K6}). The analysis of classical solutions of the system has been done mostly in \cite{K3}, where there is also a maximum principle proved for $\infty$-Harmonic maps in $2\by 2$ dimensions, namely for solutions $u:\Om \sub \R^2 \larrow \R^2$ of $\De_\infty u=0$. 

In this note we consider the question of satisfaction of the Maximum Principle for smooth solutions of the tangential part of the $\infty$-Laplace system when $N\geq2$ and $n\geq 1$
\beq \label{eq1}
Du \ot Du :D^2u\, =\, 0,\ \ \ u\ :\ \Om \sub \R^n\larrow \R^N,
\eeq
and also for smooth solutions of the single $\infty$-Laplace equation, perturbed by a linear smooth first order term
\beq \label{eq2}
Dv\ot Dv:D^2v\, +\, Dv\cdot DF\, =\, 0, \ \ \ v\ :\ \Om \sub \R^n\larrow \R.
\eeq
For the case of the system, the Maximum Principle for the vectorial solution is considered in the sense of the standard Maximum Principle for either the projections along lines, or for the modulus of the solution. We show that in both cases, the Maximum Principle in general fails. The case of projections has an elegant restatement called the Convex Hull Property and is well-known in Calculus of Variations (see e.g.\ \cite{K1}):
\[
u(\Om)\, \sub \, \co\big( u(\p\Om)\big).
\]
This inclusion says that the image is contained in the closed convex hull of the boundary values. The equivalence of this form with the form
\[
\underset{\Om}\sup\, \xi \cdot u\, \leq\, \underset{\p \Om}\max\, \xi \cdot u,\ \ \ \forall\, \xi\in \R^N,
\]
can be seen by writting the convex hull as the intersections of all affine halfspaces containing it. The case of the single equation is very interesting because is in direct constrast to the case of uniformly elliptic operators, where perturbation by a first order term does not violate the Maximum Principle.

Our examples are the following:

\begin{example}[\textbf{Failure of the Convex Hull Property for $u$ when $Du\ot Du :D^2u$ $=0$}] \label{ex1} 

Let $N\geq 2, n\geq 1$. 

\ms

\noi a) There exists $\xi \in \R^N$ and $u^1: \R^n  \larrow \R^N$ in $C^\infty(\R^n )^N$ such that $u^1$ solves the same system on $\R^n$ but the projection $\xi\cdot u^1$ along the line spanned by $\xi$ does not satisfy neither the Maximum nor the Minimum Principle:
\[
\exists \ \Om^\pm \sub \R^n \ : \ \ \ \ \begin{array}{l}
\underset{\Om^+}\sup\, \xi \cdot u^1\, \not\leq\, \underset{\p \Om^+}\max\, \xi \cdot u^1\ms\\
\underset{\Om^-}\inf\, \xi \cdot u^1\, \not\geq\, \underset{\p \Om^-}\min \,\xi \cdot u^1.
\end{array}
\]

\noi b) There exists $u^2: \R^n\set\{0\} \larrow \R^N$ in $C^\infty(\R^n\set\{0\} )^N$ which solves
\[
Du \ot Du :D^2u\, =\, 0, \ \ \text{ on }\R^n\set\{0\} ,
\]
such that $u^2$ does not satisfy the Convex Hull Property over a bounded domain:
\[
\exists\ \Om \Subset \R^n\set \{0\}\ : \ \ \ \ u^2(\Om) \not\sub \co \big(u^2(\p \Om)\big).  
\]
\end{example}

\begin{example}[\textbf{Failure of the Maximum Principle for $|u|$ when $Du\ot Du :D^2u$ $=0$}] \label{ex2}

Let $N\geq 2, n\geq 1$.  There exists $u^3: \R^n \larrow \R^N$ in $C^\infty(\R^n)^N$ which solves
\[
Du \ot Du :D^2u\, =\, 0, \ \ \text{ on }\R^n,
\]
such that $|u^3|$ does not satisfy the Maximum Principle:
\[
\exists \ \Om \sub \R^n \ : \ \ \ \ \sup_{\Om} |u^3|\, \not\leq\, \max_{\p \Om} |u^3|.
\]
\end{example}

The construction of this last example can be modified in order to give bounded domain (along the lines of the constructions that will be made for b) of Example \ref{ex1}), but we will refrain from providing all the details in this case. The same comment applies to the Example \ref{ex3} that follows.

\begin{corollary}
The above two examples show that in general it is not possible the Maximum Principle to be true for the tangential part \eqref{eq1} of $\De_\infty$ in positive codimension $N-n>0$, without utilising the \emph{complete} $\infty$-Laplace PDE system \eqref{1}.
\end{corollary}

In particular, it was shown in \cite{K2} that when $n=1\leq N$, the tangential system \eqref{eq1} is (for classical solutions) equivalent to the poperty that the solution has constant speed. On the other hand, it was shown that the full system  \eqref{1} is equivalent to that the solution is affine.

\begin{example}[\textbf{Failure of the Maximum Principle for $v$ when $Dv\ot Dv :D^2v+Dv$ $\cdot DF=0$}] \label{ex3}

Let $n\geq 1$. There exist $F\in C^\infty(\R^n)$ such that the solution $v\in C^\infty(\R^n)$ to
\[
Dv\ot Dv :D^2v\, +\, Dv \cdot DF\, =\, 0, \ \ \text{ on }\R^n,
\]
does not satisfy neither the Maximum nor the Minimum Principle:
\[
\exists \ \Om^\pm \sub \R^n \ : \ \ \ \ \begin{array}{l}
\underset{\Om^+}\sup\, v\, \not\leq\, \underset{\p \Om^+}\max\, v,\ms\\
\underset{\Om^-}\inf\, v\, \not\geq\, \underset{\p \Om^-}\min \,v.
\end{array}
\]
\end{example}

By taking $n=1$ in the above examples, all the domains that will be exhibited become by construction bounded.

\begin{remark} We note that in the above examples the non-homogenous single $\infty$-Laplace equation \eqref{eq2} and the tangential $\infty$-Laplace system \eqref{eq1} are closely related and this fact is heavily used in our constructions: for any vector solution $u$ of \eqref{eq1}, after a calculation it can be seen that both the projection $\xi \cdot u$ along a fixed direction $\xi$ and also the modulus $|u|$ satisfy a non-homogeneous equation of the form \eqref{eq2}.
\end{remark}

We conclude this introduction by noting that very recently in \cite{K8} the first author introduced a new theory of generalised solutions and proved existence to the Dirichet problem for the $\infty$-Laplacian. See also \cite{K9,K10}.

\section{Constructions}

\noi \textbf{Example \ref{ex1}.} a) Let $w_1 \in C^\infty_c(\R)$ be the function
\[
w_1(t)\, :=\, \left\{
\begin{array}{l}
-\, e^{\dfrac{1}{(1-t)^2-1}},\ \ \ \ \ 0<t<2,\ms\\

+\, e^{\dfrac{1}{(1+t)^2-1}},\ \ \ -2<t<0,\ms\ms\\

0, \hspace{75pt} \text{ otherwise}.
\end{array}
\right.
\] 

\[
\underset{\text{Figure 1. Illustration of the function $w_1$. }}{\includegraphics[scale=0.2]{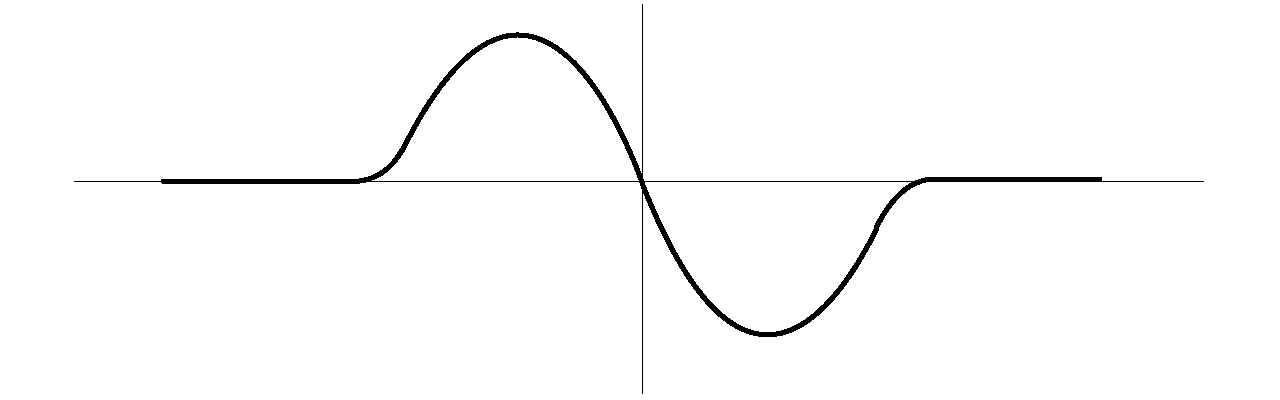}}
\]

Fix also 
\[
M\, >\, \sup_{\R}|w_1'| 
\]
and define $w_2\in C^\infty(\R)$ by taking
\[
w_2(t)\, :=\, \int_0^t \sqrt{M^2\, -\, |w_1'(s)|^2}\, ds.
\]
Then, the map $\ga: \R \larrow \R^2$ given by
\[
\ga(t)\, :=\, \big(w_1(t), w_2(t)\big)^\top
\]
is a constant speed curve and in $C^\infty(\R)^2$:
\[
|\ga'(t)|^2\, =\, |w_1'(t)|^2\, +\, |w_2'(t)|^2\, =\, M^2, \ \ \ t\in \R.
\]
Then, we define $u^1: \R^n \larrow \R^2$ by
\[
u^1(x)\, :=\, \ga(x_1), \ \ \ \ x=(x_1,...,x_n)^\top\in \R^n,\ n\geq1.
\]
The smooth map $u^1$ is a classical global solution of the tangential part of the $\infty$-Laplacian: since 
\[
Du^1(x)\, =\, \big(\ga'(x_1),0,...,0 \big)^\top,\ \ \ \ |Du^1(x)|^2\,=\, |\ga'(x_1)|^2,
\]
we have
\begin{align}
Du^1(x) \ot Du^1 (x):D^2u^1(x)\, &=\, Du^1(x)\, D\left(\frac{1}{2}|Du^1|^2 \right)(x) \nonumber\\
&=\, \ga'(x_1)\, \left(\frac{1}{2}|\ga'|^2 \right)'(x_1) \nonumber\\
&=\, 0, \nonumber
\end{align}
for all $x\in \R^n$. However, for
\[
\Om^-\,=\, (-2,0)\by \R, \ \ \ \Om^+\,=\, (0,2)\by \R,
\]
and for 
\[
\xi\, =\, e^1\, =\, (1,0)^\top, 
\]
we have
\[
\sup_{\Om^+} \xi\cdot u^1 \, >\,  \max_{\p \Om^+}\xi\cdot u^1 \,=\, 0\, =\, \min_{\p \Om^-}\xi\cdot u^1 \, >\, \inf_{\Om^-} \xi\cdot u^1 .
\]
Hence, $\xi \cdot u^1$ does not satisfy neither the maximum nor the minimum principle. By writing a convex set as the intersection of all halfspaces containing it, it follows that $u^1$ can not satisfy the convex hull property either.
\[
\underset{\text{Figure 2. Illustration of the Convex Hull Property. }}{\includegraphics[scale=0.2]{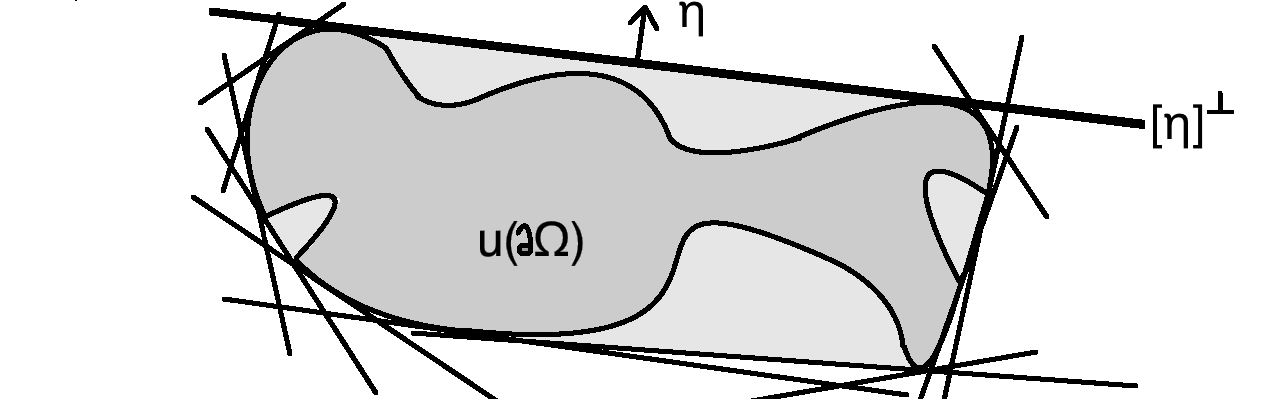}}
\]
\ms

b) We can also modify the previous construction of a) in order to have bounded domains $\Om^\pm$. Indeed, define  
\[
z_1(t)\, :=\, \left\{
\begin{array}{l}
e^{\dfrac{1}{(2-t)^2-1}},\ \ \ \ \ \ 1<t<3,\ms\\

0, \hspace{70pt} \text{ otherwise},
\end{array}
\right.
\] 
and choose as before
\[
M\, >\, \sup_{\R}|z_1'| ,
\]
\[
z_2(t)\, :=\, \int_0^t \sqrt{M^2\, -\, |z_1'(s)|^2}\, ds.
\]
Then, we define 
\[
u^2\ :\ \R^n\set\{0\} \larrow \R^2
\]
in $C^\infty(\R^n\set\{0\})^2$ by taking
\[
u^2(x)\, :=\, \big(z_1(|x|),z_2(|x|)\big)^\top,\ \ \ x\in \R^n\set\{0\},\ n\geq 1.
\]
Then, $u^2$ is a smooth solution of the vectorial Eikonal equation: indeed, since 
\[
Du^2(x)\, =\, \big(z'_1(|x|),z'_2(|x|))^\top \ot \frac{x}{|x|}, \ \ \ x\neq 0,
\]
we have
\begin{align}
\big|Du^2(x)\big|^2\, &=\, \Big|z_1'(|x|)\frac{x}{|x|} \Big|^2\, +\, \Big|z_2'(|x|)\frac{x}{|x|} \Big|^2 \nonumber\\
 &=\, \big|z_1'(|x|) \big|^2\, +\, \big|z_2'(|x|)\big|^2 \nonumber\\
 &=\, \big|z_1'(|x|) \big|^2\, +\, \big(M^2\, -\, \big|z_1'(|x|)\big|^2\big) \nonumber\\
&=\, M^2, \nonumber
\end{align}
for any $x\in \R^n\set\{0\}$. Hence, $u^2$ is a classical solution of the tangential part of the $\infty$-Laplacian:
\[
Du^2 \ot Du^2 :D^2u^2\, =\, Du^2\, D\left(\frac{1}{2}\big|Du^2\big|^2\right)\,=\, 0.
\]
However, for
\[
\Om\, :=\, \{x\in \R^n\ :\ 1<|x|<3\}
\]
and 
\[
\xi \, :=\, e^1\, =\, (1,0)^\top,
\]
we have
\[
\sup_{\Om} \xi\cdot u^2 \, >\,  \max_{\p \Om}\xi\cdot u^2 \,=\, 0.
\]
Hence, the Convex Hull Property fails for $u^2$ on $\Om$.

\ms
\ms

\noi \textbf{Example \ref{ex2}.} Suppose $u:\R^n \larrow \R^N$ is in $C^\infty(\R^n)^N$ and satisfies $|u|>0$. Then, we write in polar coordinates 
\[
u\, =\, \rho \, n,\ \ \ \rho\,:=\, |u|,\ \ n\,:=\,\frac{u}{|u|},
\]
and calculate
\[
Du\, =\, n\ot D\rho\, +\, \rho Dn.
\]
Since $|n|^2=1$, we have by differentation $n^\top Dn=0$, that is $n \, \bot\, D_i n$, $i=1,...,n$. Hence, by perpendicularity
\[
|Du|^2\, =\, |D\rho|^2\, +\, \rho^2|Dn|^2.
\]
Then, we may also calculate
\begin{align}
Du\ot Du:D^2u\, &=\, Du\, D\left(\frac{1}{2}|Du|^2\right) \nonumber\\
 &=\, \frac{1}{2}\Big(n\ot D\rho\, +\, \rho Dn \Big)D\Big(|D\rho|^2\, +\, \rho^2|Dn|^2 \Big).\nonumber
\end{align}
We now define functions 
\[
\rho^* \in C^\infty(\R^n),\ \ \ n^*\in C^\infty(\R^n)^N
\]
with
\[
\rho^*\, >\, 0 \ \text{ on }\R^n,\ \ \ |n^*|^2\, \equiv \,1\ \text{ on }\R^n,
\]
and show that if we define $u^3$ by taking $\rho^*$ and $n^*$ as its polar and radial coordinates, then $u^3$ is a global solution of \eqref{eq1} on $\R^n$, but $|u^3|$ fails to satisfy the Maximum Principle. Indeed, we define 
\[
\rho^*(t)\, :=\, e^{-|t|^2}
\]
\[
\underset{\text{Figure 3. Illustration of the function $\rho^*$. }}{ \includegraphics[scale=0.2]{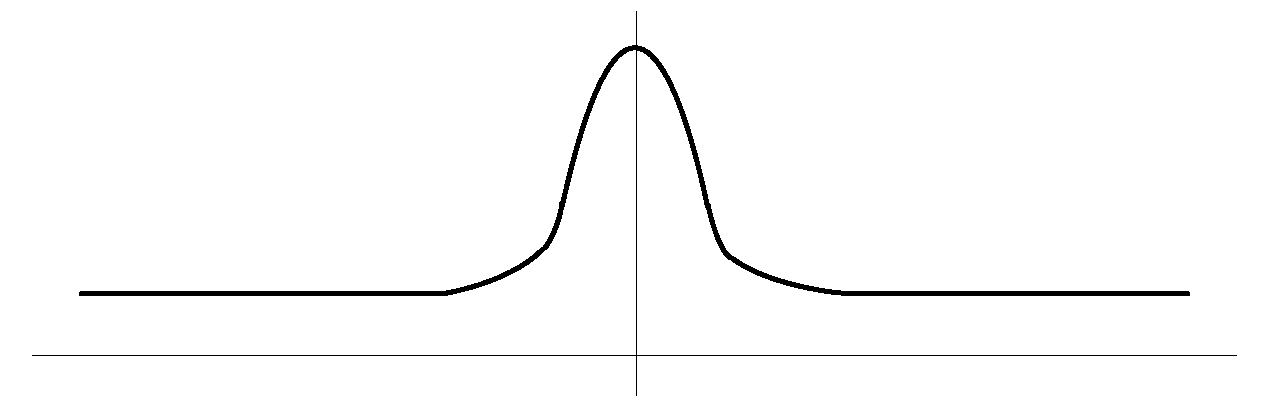}}
\]
and fix a constant
\[
M\, >\, \sup_{\R}|(\rho^*)'|.
\]
Let now $\ga : \R \larrow \R^2$ be a smooth unit speed curve parameterising the circle, e.g.\ we can take
\[
\ga(t)\, :=\, (\cos t, \sin t)^\top.
\]
Then, we define
\[
K(t):= \int_0^t\sqrt{\frac{M^2\, -\, \big| (\rho^*)'(s) \big|^2}{(\rho^*)^2(s)}}\,ds
\]
and
\[
n^*(t)\, :=\, \ga\big(K(t)\big).
\]
The previous definitions imply that $\rho^*,n^*$ satisfy
\[
(n^*)'(t)\, =\, \ga'\big(K(t)\big)\, K'(t)
\]
and hence
\[
|(n^*)'(t)|\,=\, \big| \ga'\big(K(t)\big)\, K'(t) \big|\,=\, \big| K'(t) \big|,
\]
while,
\[
 \big| K'(t) \big|^2\, =\, \frac{M^2\, -\, \big| (\rho^*)'(t) \big|^2}{(\rho^*)^2(t)}
\]
which (by combining the previous two equalities) implies
\[
(\rho^*)^2(t) \big| (n^*)'(t) \big|^2\, +\, \big|(\rho^*)'(t) \big|^2\, =\, M^2,
 \]
for all $t\in \R$. Hence, if 
\[
u^3\ :\ \R^n \larrow \R^2
\]
is given by 
\[
u^3(x)\, :=\, \rho^*(x_1) n^*(x_1), \ \ \ x=(x_1,...,x_n)^\top\in \R^n,
\]

\ms
\noi we have that $u^3$ is smooth and since $|Du^3|\equiv M$, it solves 

\[
Du^3\ot Du^3:D^2u^3\, =\, 0,\ \ \text{ on }\R^n.
\]
However, for
\[
\Om\, :=\, \big\{ x\in \R^n\ : \ |x_1|<1 \big\},
\]
we have that $|u^3|=\rho^*$ does not satisfy the maximum principle:
\[
\sup_{\Om} |u^3|\, =\,1\, >\, \frac{1}{e}\,=\, \max_{\p\Om}|u^3|.
\]

\ms
\ms

\noi \textbf{Example \ref{ex3}.} In the setting of Example \ref{ex1}, consider again the function $w_1 \in C^\infty_c(\R)$ given by
\[
w_1(t)\, :=\, \left\{
\begin{array}{l}
-\, e^{\dfrac{1}{(1-t)^2-1}},\ \ \ \ \ 0<t<2,\ms\\

+\, e^{\dfrac{1}{(1+t)^2-1}},\ \ \ -2<t<0,\ms\ms\\

0, \hspace{75pt} \text{ otherwise}.
\end{array}
\right.
\] 
Fix 
\[
M\, >\, \sup_{\R}|w_1'| 
\]
and define $w_2\in C^\infty(\R)$ by taking
\[
w_2(t)\, :=\, \int_0^t \sqrt{M^2\, -\, |w_1'(s)|^2}\, ds.
\]
Then, the map $\ga: \R \larrow \R^2$ given by
\[
\ga(t)\, :=\, \big(w_1(t), w_2(t)\big)^\top
\]
is a constant speed curve and in $C^\infty(\R)^2$:
\[
|\ga'(t)|^2\, =\, |w_1'(t)|^2\, +\, |w_2'(t)|^2\, =\, M^2, \ \ \ t\in \R.
\]
Then, we define $v: \R^n \larrow \R$ in $C^\infty(\R^n)$ by
\[
v(x)\, :=\, w_1(x_1), \ \ \ \ x=(x_1,...,x_n)^\top\in \R^n,\ n\geq1.
\]
We also define $F: \R^n \larrow \R$ in $C^\infty(\R^n)$ by
\[
F(x)\, :=\, \frac{1}{2}\big|w'_2(x_1)\big|^2, \ \ \ \ x=(x_1,...,x_n)\in \R^n,\ n\geq1.
\]
Given $v,F$ as above, we calculate
\begin{align}
Dv \ot Dv :D^2v\, +\, Dv\cdot DF\, &=\, w_1'w_1'w_1''\ +\ w_1' \left( \frac{1}{2}\big|w'_2\big|^2 \right)' \nonumber\\
&=\, \frac{1}{2} w_1'\left( \big|w'_1\big|^2 \right)'\, +\, \frac{1}{2}w_1' \left( M^2\, -\, \big|w'_1\big|^2 \right)'  \nonumber\\
&=\, 0, \nonumber
\end{align}
on $\R^n$. Hence, $v$ is a global smooth solution of
\[
Dv \ot Dv :D^2v\, +\, Dv\cdot DF\, =\, 0,
\]
on $\R^n$. However, for
\[
\Om^-\,:=\, (-2,0)\by \R, \ \ \ \Om^+\,:=\, (0,2)\by \R,
\]
we have
\[
\sup_{\Om^+} \, v \, >\,  \max_{\p \Om^+}\, v \,=\, 0\, =\, \min_{\p \Om^-}\, v \, >\, \inf_{\Om^-} \, v .
\]
Hence, $v$ does not satisfy neither the Maximum nor the Minimum Principle.

\ms\ms

\noi \textbf{Acknowledgement.} The second author has been financially partially supported by 
NSF award DMS-1001179.
\ms

\bibliographystyle{amsplain}

\end{document}